\documentclass{siamltex}
\usepackage{amsmath,amsfonts,amssymb}
\usepackage[usenames]{color}

\usepackage[mathscr]{eucal}

\usepackage{stmaryrd}

\usepackage{amsbsy}

\usepackage{bm}

\usepackage{xspace}

\usepackage{paralist}

\usepackage{ifpdf}

\usepackage{graphicx}
\ifpdf
\usepackage{epstopdf}
\else
\usepackage{epsfig}
\fi

\usepackage{braket}

\usepackage{subfig}

\ifpdf
\usepackage[colorlinks,urlcolor=blue,citecolor=blue,linkcolor=blue]{hyperref}
\else
\newcommand{\href}[2]{{#2}}
\usepackage{url}
\fi

\usepackage{tabulary}

\usepackage{algorithm}
\usepackage{algpseudocode}

\ifpdf
\newcommand{\Sec}[1]{\hyperref[sec:#1]{\S\ref*{sec:#1}}} 
\newcommand{\App}[1]{\hyperref[sec:#1]{Appendix~\ref*{sec:#1}}} 
\newcommand{\Eqn}[1]{\hyperref[eq:#1]{{\rm (\ref*{eq:#1})}}} 
\newcommand{\Part}[1]{\hyperref[part:#1]{(\ref*{part:#1})}} 
\newcommand{\Fig}[1]{\hyperref[fig:#1]{Figure~\ref*{fig:#1}}} 
\newcommand{\Tab}[1]{\hyperref[tab:#1]{Table~\ref*{tab:#1}}} 
\newcommand{\Thm}[1]{\hyperref[thm:#1]{Theorem~\ref*{thm:#1}}} 
\newcommand{\Lem}[1]{\hyperref[lem:#1]{Lemma~\ref*{lem:#1}}} 
\newcommand{\Prop}[1]{\hyperref[prop:#1]{Property~\ref*{prop:#1}}} 
\newcommand{\Cor}[1]{\hyperref[cor:#1]{Corollary~\ref*{cor:#1}}} 
\newcommand{\Def}[1]{\hyperref[def:#1]{Definition~\ref*{def:#1}}} 
\newcommand{\Alg}[1]{\hyperref[alg:#1]{Algorithm~\ref*{alg:#1}}} 
\newcommand{\Ex}[1]{\hyperref[ex:#1]{Example~\ref*{ex:#1}}} 
\newcommand{\As}[1]{\hyperref[as:#1]{Assumption~{\rm\ref*{as:#1}}}} 
\newcommand{\Reg}[1]{\hyperref[as:#1]{Condition~\ref*{reg:#1}}} 
\newcommand{\AlgLine}[2]{\hyperref[alg:#1]{line~\ref*{line:#2} of Algorithm~\ref*{alg:#1}}}
\newcommand{\AlgLines}[3]{\hyperref[alg:#1]{lines~\ref*{line:#2}--\ref*{line:#3} of Algorithm~\ref*{alg:#1}}}
\else
\newcommand{\Sec}[1]{{\S\ref{sec:#1}}} 
\newcommand{\App}[1]{{Appendix~\ref{sec:#1}}} 
\newcommand{\Eqn}[1]{{(\ref{eq:#1})}} 
\newcommand{\Part}[1]{{(\ref{part:#1})}} 
\newcommand{\Fig}[1]{{Figure~\ref{fig:#1}}} 
\newcommand{\Tab}[1]{{Table~\ref{tab:#1}}} 
\newcommand{\Thm}[1]{{Theorem~\ref{thm:#1}}} 
\newcommand{\Lem}[1]{{Lemma~\ref{lem:#1}}} 
\newcommand{\Prop}[1]{{Property~\ref{prop:#1}}} 
\newcommand{\Cor}[1]{{Corollary~\ref{cor:#1}}} 
\newcommand{\Def}[1]{{Definition~\ref{def:#1}}} 
\newcommand{\Alg}[1]{{Algorithm~\ref{alg:#1}}} 
\newcommand{\Ex}[1]{{Example~\ref{ex:#1}}} 
\newcommand{\Reg}[1]{{R~\ref*{reg:#1}}} 
\fi

\usepackage{sudoku}
\usepackage{booktabs}

\usepackage{tikz}
\usetikzlibrary{shapes,snakes}
\usepackage{xcolor}

\usepackage{url}
\usepackage{multirow}

\newcommand{\Real}{\mathbb{R}}

\newcommand{\Tra}{^{\sf T}} 

\newcommand{\V}[1]{{\bm{\mathbf{\MakeLowercase{#1}}}}} 
\newcommand{\VE}[2]{\MakeLowercase{#1}_{#2}} 

\newcommand{\M}[1]{{\bm{\mathbf{\MakeUppercase{#1}}}}} 
\newcommand{\Mhat}[1]{{\bm{\hat \mathbf{\MakeUppercase{#1}}}}} 

\newcommand{\T}[1]{\boldsymbol{\mathscr{\MakeUppercase{#1}}}} 
\newcommand{\TE}[2]{\MakeLowercase{#1}_{#2}} 

\begin{document}
\title{Techniques for Solving Sudoku Puzzles}
\author{Eric C. Chi\footnotemark[1] \and Kenneth Lange\footnotemark[1] \footnotemark[2]}
\maketitle

\renewcommand{\thefootnote}{\fnsymbol{footnote}}
\footnotetext[1]{Dept. Human Genetics, University of California, Los Angeles, CA. Email: ecchi@ucla.edu}
\footnotetext[2]{Depts. of Biomathematics and Statistics, University of California, Los Angeles,
CA. Email: klange@ucla.edu}
\renewcommand{\thefootnote}{\arabic{footnote}}

\begin{abstract}
Solving Sudoku puzzles is one of the most popular pastimes in the world.  Puzzles range in difficulty from easy to very challenging; the hardest puzzles tend to have the most empty cells.  The current paper explains and compares three algorithms for solving Sudoku puzzles. Backtracking, simulated annealing, and alternating projections are generic methods for attacking combinatorial optimization problems.  Our results
favor backtracking.  It infallibly solves a Sudoku puzzle or deduces that a unique solution does not exist. However, backtracking does not scale well in high-dimensional combinatorial optimization.  Hence, it is useful to expose students in the mathematical sciences to the other two solution techniques in a concrete setting.  Simulated annealing shares a common structure with MCMC (Markov chain Monte Carlo) and enjoys wide applicability.  The method of alternating projections solves the feasibility problem in  convex programming. Converting a discrete optimization problem into a continuous optimization problem opens up the possibility of handling combinatorial problems of much higher dimensionality.
\end{abstract}


\section{Introduction}

As all good mathematical scientists know, a broad community has contributed to the invention of modern algorithms. Computer scientists, applied mathematicians, statisticians, economists, and physicists, to name just a few, have made lasting contributions. Exposing students to a variety of perspectives outside the realm of their own disciplines sharpens their instincts for modeling and arms them with invaluable tools.  In this spirit, the current paper discusses techniques for solving Sudoku puzzles, one of the most popular pastimes in the world.  One could make the same points with more serious applications, but it is hard to imagine a more beguiling introduction to the algorithms featured here. Sudoku diagrams are special cases of the Latin squares long familiar in experimental design and, as such, enjoy interesting mathematical and statistical properties \cite{BaiCamCon2008}. The complicated constraints encountered in solving Sudoku puzzles have elicited many clever heuristics that amateurs use to good effect. Here we examine three generic methods with broader scientific and societal applications.  The fact that one of these methods outperforms the other two is mostly irrelevant.  No two problem categories are completely alike, and it is best to try many techniques before declaring a winner.

The three algorithms tested here are simulated annealing, alternating projections, and backtracking.   Simulating annealing is perhaps the most familiar to statisticians.  It is the optimization analog of MCMC (Markov chain Monte Carlo) and has been employed to solve a host of combinatorial problems.  The method of alternating projections was first proposed by von Neumann \cite{Neu1950} to find a feasible point in the intersection of a family of hyperplanes.  Modern versions of alternating projections more generally seek a point in the intersection of a family of closed convex sets.  Backtracking is a standard technique taken from the toolkits of applied mathematics and computer science.  Backtracking infallibly finds all solutions of a Sudoku puzzle or determines that no solution exists. Its Achilles heel of excessive computational complexity does not come into play with Sudoku puzzles  because they are, despite appearances, relatively benign computationally.   Sudoku puzzles are instances of the satisfiability problem in computer science.  As problem size increases, such problems are combinatorially hard and often defy backtracking.  For this reason alone, it is useful to examine alternative strategies.

In a typical Sudoku puzzle, there are 81 cells arranged in a 9-by-9 grid, some of which are occupied by numerical clues. See Figure~\ref{fig:sample_game}. The goal is to fill in the remaining cells subject to the following three rules:
\begin{figure}
\centering
\begin{tikzpicture}[scale=0.75]
    \draw[black!50] (0,0) grid (9,9);
	\draw[black!100] (0,0) rectangle (3,3);
	\draw[black!100] (3,0) rectangle (6,3);
	\draw[black!100] (6,0) rectangle (9,3);
	\draw[black!100] (0,3) rectangle (3,6);
	\draw[black!100] (3,3) rectangle (6,6);
	\draw[black!100] (6,3) rectangle (9,6);
	\draw[black!100] (0,6) rectangle (3,9);
	\draw[black!100] (3,6) rectangle (6,9);
	\draw[black!100] (6,6) rectangle (9,9);	
	\node[regular polygon, regular polygon sides=4] at (3.5,0.5) {$5$};	
	\node[regular polygon, regular polygon sides=4] at (5.5,0.5) {$6$};
	\node[regular polygon, regular polygon sides=4] at (6.5,0.5) {$8$};
	\node[regular polygon, regular polygon sides=4] at (7.5,0.5) {$2$};	

	\node[regular polygon, regular polygon sides=4] at (3.5,1.5) {$3$};
	\node[regular polygon, regular polygon sides=4] at (5.5,1.5) {$4$};	
	\node[regular polygon, regular polygon sides=4] at (7.5,1.5) {$1$};
	\node[regular polygon, regular polygon sides=4] at (8.5,1.5) {$9$};
		
	\node[regular polygon, regular polygon sides=4] at (1.5,2.5) {$8$};
	\node[regular polygon, regular polygon sides=4] at (4.5,2.5) {$2$};	
	\node[regular polygon, regular polygon sides=4] at (5.5,2.5) {$1$};

	\node[regular polygon, regular polygon sides=4] at (4.5,3.5) {$7$};	
	\node[regular polygon, regular polygon sides=4] at (8.5,3.5) {$1$};
	
	\node[regular polygon, regular polygon sides=4] at (0.5,4.5) {$2$};
	\node[regular polygon, regular polygon sides=4] at (3.5,4.5) {$8$};	
	\node[regular polygon, regular polygon sides=4] at (4.5,4.5) {$6$};
	\node[regular polygon, regular polygon sides=4] at (7.5,4.5) {$7$};
	\node[regular polygon, regular polygon sides=4] at (8.5,4.5) {$4$};	
	
	\node[regular polygon, regular polygon sides=4] at (0.5,5.5) {$7$};
	\node[regular polygon, regular polygon sides=4] at (2.5,5.5) {$4$};
	\node[regular polygon, regular polygon sides=4] at (3.5,5.5) {$1$};	
	\node[regular polygon, regular polygon sides=4] at (5.5,5.5) {$5$};
	\node[regular polygon, regular polygon sides=4] at (6.5,5.5) {$2$};
	
	\node[regular polygon, regular polygon sides=4] at (1.5,6.5) {$3$};
	\node[regular polygon, regular polygon sides=4] at (2.5,6.5) {$2$};
	\node[regular polygon, regular polygon sides=4] at (3.5,6.5) {$9$};
	\node[regular polygon, regular polygon sides=4] at (6.5,6.5) {$1$};			
	\node[regular polygon, regular polygon sides=4] at (7.5,6.5) {$4$};

	\node[regular polygon, regular polygon sides=4] at (1.5,7.5) {$4$};

	\node[regular polygon, regular polygon sides=4] at (0.5,8.5) {$1$};
	\node[regular polygon, regular polygon sides=4] at (1.5,8.5) {$5$};
	\node[regular polygon, regular polygon sides=4] at (2.5,8.5) {$7$};
	\node[regular polygon, regular polygon sides=4] at (3.5,8.5) {$6$};	
	\node[regular polygon, regular polygon sides=4] at (4.5,8.5) {$4$};			
	\node[regular polygon, regular polygon sides=4] at (7.5,8.5) {$8$};	
\end{tikzpicture}
\caption{Sample Puzzle\label{fig:sample_game}}
\end{figure}

\begin{itemize}
\item[1.] Each integer between 1 and 9 must appear exactly once in a row,
\item[2.] Each integer between 1 and 9 must appear exactly once in a column,
\item[3.] Each integer between 1 and 9 must appear exactly once in each of the 3-by-3 subgrids.
\end{itemize}

Solving a Sudoku game is a combinatorial task of intermediate complexity. The general problem of filling in an incomplete $n^2 \times n^2$ grid with $n \times n$ subgrids belongs to the class of NP-complete problems \cite{YatSet2003}. These problems are conjectured to increase in computational complexity at an exponential rate in $n$.  Nonetheless,  a well planned exhaustive search can work quite well for a low value of $n$ such as $9$. For larger values of $n$, brute force, no matter how cleverly executed, is simply not an option. In contrast, simulated annealing and alternating projections may yield good approximate solutions and partially salvage the situation.

In the rest of this paper, we describe the three methods for solving Sudoku puzzles and compare them on a battery of puzzles. The puzzles range in difficulty from pencil and paper exercises to hard benchmark tests that often defeat the two approximate methods. Our discussion reiterates the rationale for equipping students with the best computational tools.

\section{Three methods for solving Sudoku}
\label{sec:methods}

\subsection{Backtracking}

Backtracking systematically grows a partial solution until it becomes a full solution or violates a constraint \cite{Ski2008}. In the latter case it backtracks to the next permissible partial solution and begins the growing process anew.  The advantage of backtracking is that a block of potential solutions can be discarded en masse. Backtracking starts by constructing for each empty Sudoku cell 
$(i,j)$ a list $L_{ij}$ of compatible digits.  This is done by scanning the cell's row, column, and subgrid. The empty cells are then
ordered by the cardinalities of the lists $|L_{ij}|$. For example in Figure~\ref{fig:sample_game}, two cells $(7,4)$ and $(9,5)$ possess lists $L_{74}=\{7\}$ and $L_{95}=\{9\}$ with cardinality 1 and come first.  Next come cells such as $(1,6)$ with $L_{16}=\{2,3\}$, $(1,7)$ with $L_{17}=\{3,9\}$, and $(1,9)$ with $L_{19}=\{2,3\}$ whose lists have cardinality 2.  Finally come cells such as $(2,9)$ with $L_{29}=\{2,3,5,6,7\}$ whose lists have maximum cardinality 5.  Partial solutions are character strings such as $s_{74}s_{95}s_{16}s_{17}s_{19}$ taken in dictionary order with the alphabet at cell $(i,j)$ limited to the list $L_{ij}$.  In dictionary order a string such as $7939$ is treated as coming after a string such as $79232$.

Backtracking starts with the string $7$ by taking the only element of $L_{74}$, grows it to $79$ by taking the only element of $L_{95}$, grows it to $792$ by taking the first element of $L_{16}$, grows it to $7923$ by taking the first element of $L_{17}$, and finally grows it to $79232$ by taking the first element of $L_{19}$. At this juncture a row violation occurs, namely a 2 in both cells $(1,6)$ and $(1,9)$. Backtracking discards all strings beginning with $79232$ and moves on to the string $79233$ by replacing the first element of $L_{19}$ by the second element of $L_{19}$.  This leads to another row violation with a 3 in both cells $(1,7)$ and $(1,9)$.  Backtracking moves back to the string $7929$ by discarding the fifth character of $79239$ and replacing the first element of $L_{17}$ by its second element.  This sets the stage for another round of growing.

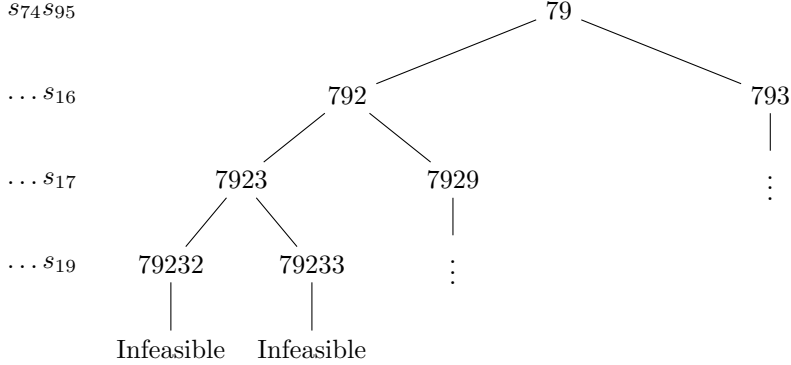
\begin{figure}
\centering
\begin{tikzpicture}[level/.style={sibling distance=75mm/#1}, scale=0.75]
\node  (z){79}
 child {node (a) {792}
  child {node  (b) {7923}
  	child {node(c) {79232}
		child{node(d) {Infeasible}
			child [grow=left] {node (q) {\quad\quad\quad\quad} edge from parent[draw=none]
			 child [grow=up] {node (r) {$\ldots\VE{s}{19}$ \quad\quad\quad\quad} edge from parent[draw=none]
			 child [grow=up] {node (s) {$\ldots\VE{s}{17}$ \quad\quad\quad\quad} edge from parent[draw=none]
			 child [grow=up] {node (t) {$\dots\VE{s}{16}$ \quad\quad\quad\quad} edge from parent[draw=none]
			 child [grow=up] {node (t) {$\VE{s}{74}\VE{s}{95}$ \quad\quad\quad\quad} edge from parent[draw=none]
			 }}}}}
		}
	}
	   child {node (e)  {79233}
	   	child {node (f) {Infeasible}
	   }}
 }
   child {node (g)  {7929}
   	child {node (h) {$\vdots$}}}
 }
  child {node (i) {793}
  child {node (j) {$\vdots$}}
  }
  ;
\end{tikzpicture}
\caption{Backtracking on the puzzle shown in Figure~\ref{fig:sample_game}. Starting from $\VE{s}{74}\VE{s}{95} = 79$, the algorithm
attempts and fails to grow the solution beyond $\VE{s}{74}\VE{s}{95}\VE{s}{16}\VE{s}{17}\VE{s}{19} = 79232$. After failing to grow the solution
beyond $\VE{s}{74}\VE{s}{95}\VE{s}{16}\VE{s}{17}\VE{s}{19} = 79233$, all partial solutions beginning with 7923 are eliminated from further consideration.
The algorithm starts anew by attempting to grow $\VE{s}{74}\VE{s}{95}\VE{s}{16}\VE{s}{17} = 7929$. }
\label{fig:backtracking}
\end{figure}

Backtracking is also known as depth first search. In this setting the strings are viewed as nodes of a tree as depicted in Figure~\ref{fig:backtracking}.  Generating strings in dictionary order constitutes a tree traversal that systematically eliminates subtrees and moves down and backs up along branches. 
Because pruning large subtrees is more efficient than pruning small subtrees, ordering of cells by cardinality compels the decision tree to have fewer branches at the top. We use the C code from Skiena and Revilla \cite{SkiRev2003} implementing backtracking on Sudoku puzzles. Backtracking has the virtue of finding all solutions when multiple solutions exist.  Thus, it provides a mechanism for validating the correctness of puzzles.

\subsection{Simulated Annealing}

Simulated annealing \cite{Cer1985, KirGelVec1983, PreTeuVet2007} attacks a combinatorial optimization problem by defining a state space of possible solutions, a cost function quantifying departures from the solution ideal, and a positive temperature parameter. For a satisfiability problem, it is sensible to equate cost to the number of constraint violations.  Solutions then correspond to states of zero cost. Each step of annealing operates by proposing a move to a new randomly chosen state.  Proposals are Markovian in the sense that they depend only on the current state of the process, not on its past history. Proposed steps that decrease cost are always accepted. Proposed steps that increase cost are taken with high probability in the early stages of annealing when temperature is high and with low probability in the late stages of annealing when temperature is low. Inspired by models from statistical physics, simulated annealing is designed to sample the state space broadly before settling down at a local minimum of the cost function.

For the Sudoku problem, a state is a $9 \times 9$ matrix (board) of integers drawn from the set $\{1,\ldots,9\}$.  Each integer appears nine times, and all numerical clues are respected. Annealing starts from any feasible board. The proposal stage randomly selects two different cells without clues. The corresponding move swaps the contents of the cells, thus preserving all digit counts. To ensure that the most troublesome cells are more likely to be chosen for swapping, we select cells non-uniformly with probability proportional to 
$\exp(i)$ for a cell involved in $i$ constraint violations.  Let $\M{b}$ denote a typical board, $c(\M{b})$ its associated cost, and 
$n$ the current iteration index. At temperature $\tau$, we decide whether to accept a proposed neighboring board $\M{b}$ by drawing a random deviate $U$ uniformly from $[0,1]$. If $U$ satisfies 
\begin{eqnarray*}
U & \le & \min\left\{\exp (\left[c(\M{b}_n)-c(\M{b})\right]/\tau_n),1\right\},
\end{eqnarray*}
then we accept the proposed move and set $\M{b}_{n+1}=\M{b}$. Otherwise, we reject the move and set $\M{b}_{n+1}=\M{b}_n$. Thus, the greater the increase in the number of constraint violations, the less likely the move is made to a proposed state. Also, the higher the temperature, the more likely a move is made to an unfavorable state.  The final ingredient of simulated annealing is the cooling schedule. In general, the temperature parameter $\tau$ starts high and slowly declines to 0, where only favorable or cost neutral moves are taken.
Typically temperature is lowered at a slow geometric rate.

\subsection{Alternating Projections}

The method of alternating projections relies on projection operators.   In the projection problem, one seeks the closest point 
$\V{x}$ in a set $C \subset \Real^d$ to a point  $\V{y} \in  \Real^d$.  Distance is quantified by the usual Euclidean norm
$\lVert \V{x} - \V{y} \rVert.$  If $\V{y}$ already lies in $C$, then the problem is trivially solved by setting $\V{x} = \V{y}$.  
It is well known that a unique minimizer exists whenever the set $C$ is closed and convex \cite{Lan2004}.  We will denote the
projection operator taking $\V{y}$ to $\V{x}$ by $P_C(\V{y}) = \V{x}$. 

\begin{figure}[t]
\centering
\includegraphics[scale=0.6]{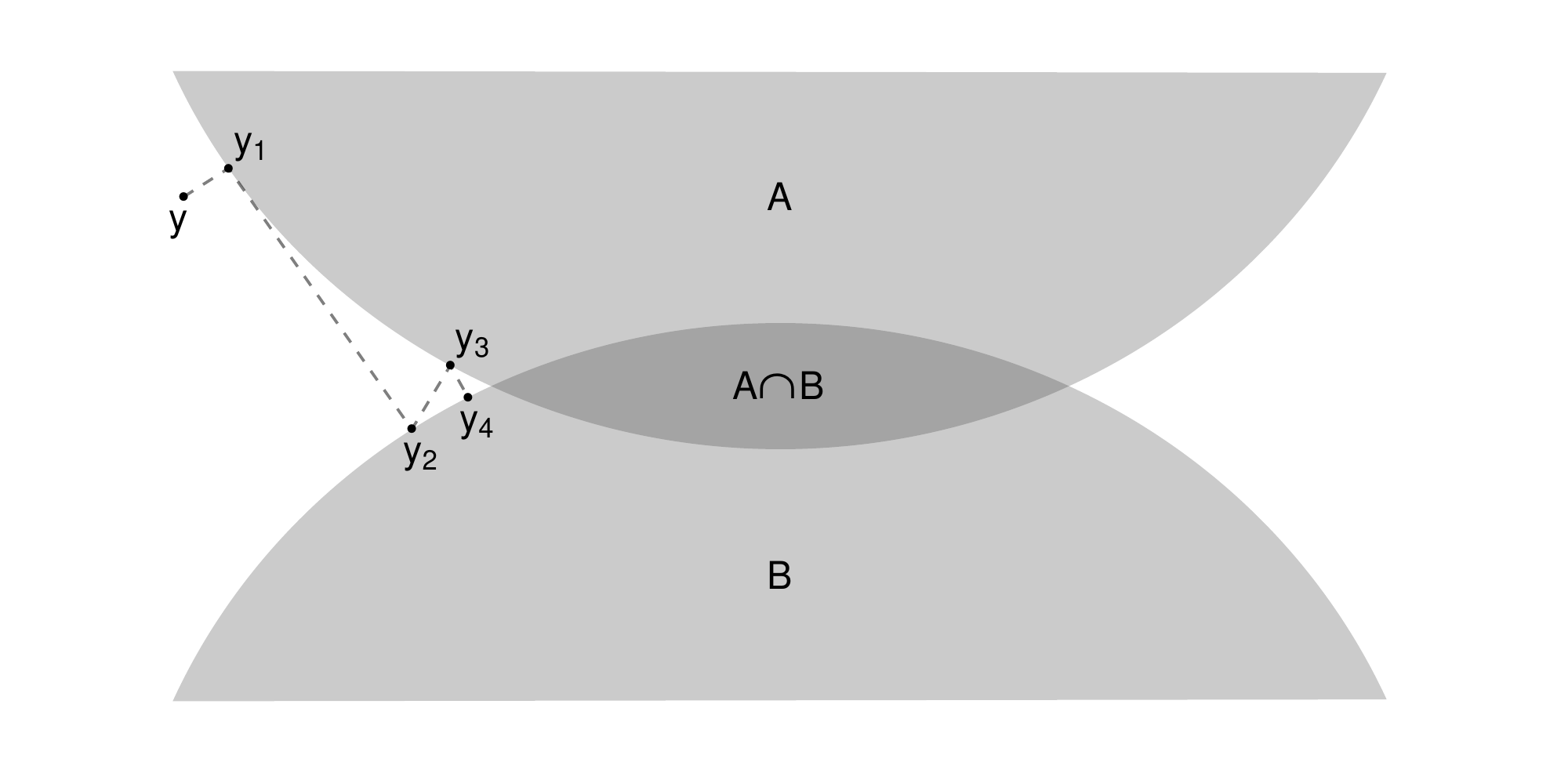}
\caption{Alternating projections find a point in $A \cap B$, where $A$ and $B$ are closed convex sets. The initial point is $\V{y}$. The sequence of points $\V{y}_n$ is generated by alternating projection onto $A$ with projection onto $B$.}
\label{fig:alternating_projections}
\end{figure}

Given a finite collection of closed convex sets with a nonempty intersection, the alternating projection algorithm finds a point in that intersection. Consider the case of two closed convex sets $A$ and $B$.  The method recursively generates a sequence $\V{y}_n$
by taking $\V{y}_0 = \V{y}$ and $\V{y}_{n+1} = P_A(\V{y}_{n})$ for $n$ even and $\V{y}_{n+1} = P_B(\V{y}_{n})$ for $n$ odd. Figure~\ref{fig:alternating_projections} illustrates a few iterations of the algorithm. As suggested by the picture, the algorithm does indeed
converge to a point in $A \cap B$ \cite{CheGol1959}.  For more than two closed convex sets with nonempty intersection, the method of alternating projections cycles through the projections in some fixed order.  Convergence occurs in this more general case as well based
on some simple theory involving paracontractive operators \cite{ElsKolNeu1992}.  The limit is not guaranteed to be the closest point in the intersection to the original point $\V{y}$.  The related but more complicated procedure known as Dykstra's algorithm \cite{Dyk1983} finds this point. 

It is easy to construct some basic projection operators.  For instance, projection onto the rectangle 
$R =  \{\V{x} \in \Real^d : a_i \le x_i \le b_i \: \mbox{for all} \: i\}$ is achieved by defining $\V{x}=P_R(\V{y})$ to have components $x_i=\min\{\max\{a_i,y_i\},b_i\}$.  This example illustrates a more general rule; namely, if $A$ and $B$ are two closed convex sets, then projection onto the Cartesian product $A \times B$ is effected by the Cartesian product operator 
$(\V{x},\V{y}) \mapsto [P_A(\V{x}),P_B(\V{y})]$. When $A$ is an entire Euclidean space, $P_A(\V{x})$ is just the identity map. Projection onto the hyperplane 
\begin{eqnarray*}
H  & = & \{\V{y} \in \Real^d : \V{v}\Tra\V{y} = c\}
\end{eqnarray*}
is implemented by the operator
\begin{eqnarray*}
P_H(\V{x}) & = & \V{x} - \frac{\V{v}\Tra\V{x} - c}{\lVert \V{v} \rVert^2} \V{v} .
\end{eqnarray*}
Projection onto the unit simplex $U = \left\{\V{x} \in \Real^d: \sum_{i=1}^d x_i = 1, \; x_i \ge 0 \:  \forall i\right\}$ is more subtle.  Fortunately there exist fast algorithms for this purpose \cite{DucShaSin2008, Mic1986}.

In either the alternating projection algorithm or Dykstra's algorithm, it is advantageous to reduce the number of participating convex sets to the minimum possible consistent with fast projection.  For instance, it is better to take the unit simplex $U$ as a whole rather than as an intersection of the halfspaces $\{\V{x}: x_i \ge 0\}$ and the affine subspace $\{\V{x}: \sum_{i=1}^d x_i =1\}$. Because our alternating projection algorithm for solving Sudoku puzzles relies on projecting onto several simplexes, it is instructive to derive the Duchi et al \cite{DucShaSin2008} projection algorithm.  Consider minimization of a convex smooth function $f(\V{x})$
over $U$.  The Karush-Kuhn-Tucker stationarity condition involves setting the gradient of the Lagrangian
\begin{eqnarray*}
{\cal L}(\V{x},\lambda,\V{\mu})  & = &  f(\V{x})+ \lambda \Big(\sum_{i=1}^d x_i-1\Big) - \sum_{i=1}^d \mu_i x_i
\end{eqnarray*}
equal to $\V{0}$.  This is stated in components as the Gibbs criterion
\begin{eqnarray*}
0  & = &  \frac{\partial}{\partial x_i} f(\V{x})+ \lambda - \mu_i 
\end{eqnarray*}
for multipliers $\mu_i \ge 0$ obeying the complementary slackness conditions $\mu_i x_i = 0$.   For the choice
$f(\V{x}) = \frac{1}{2}\|\V{x}-\V{y}\|^2$, the Gibbs condition can be solved in the form
\begin{eqnarray*}
\VE{x}{i}  & = &  \begin{cases} \VE{y}{i}-\lambda & \VE{x}{i} > 0 \\ \VE{y}{i}-\lambda+\mu_i  & \VE{x}{i} = 0. \end{cases}
\end{eqnarray*}
If we let $I_+ = \{i: \VE{x}{i}>0\}$, then the equality constraint
\begin{eqnarray*}
1 & = &   \sum_{i \in I_+} \VE{x}{i}  \:\;\, = \:\;\,  \sum_{i \in I_+} \VE{y}{i} - |I_+| \lambda 
\end{eqnarray*}
implies
\begin{eqnarray*}
\lambda & = & \frac{1}{|I_+|} \Big(\sum_{i \in I_+} \VE{y}{i} - 1 \Big) .
\end{eqnarray*}
The catch, of course, is that we do not know $I_+$.

The key to avoid searching over all $2^d$ subsets is the simple observation that the $\VE{x}{i}$ and $\VE{y}{i}$ are consistently ordered.
Suppose on the contrary that $\VE{y}{i}< y_j$ and $x_j < \VE{x}{i}$.  For small $s>0$ substitute $\VE{x}{j}+s$ for $x_j$ and $\VE{x}{i}-s$
for $\VE{x}{i}$.  The objective function $f(\V{x}) = \frac{1}{2}\lVert \V{x}-\V{y} \rVert^2$ then changes by the amount 
\begin{equation*}
 \frac{1}{2} \Big[(\VE{x}{i}-s-\VE{y}{i})^2+(x_j+s-y_j)^2-(\VE{x}{i}-\VE{y}{i})^2-(x_j-y_j)^2 \Big] = s(\VE{y}{i}-y_j+x_j-\VE{x}{i})+s^2 ,
\end{equation*}
which is negative for $s$ small. Let $\VE{w}{i}$ denote the $i$th largest entry of $\V{y}$. Then the Gibbs condition implies that
$\VE{w}{1} \geq \VE{w}{2} \geq \ldots \geq \VE{w}{\lvert I_+ \rvert} > \lambda$. Thus, to determine $\lambda$ we seek the largest $k$
such that
\begin{equation*}
	\VE{w}{k} > \frac{1}{k} \left ( \sum_{i=1}^k \VE{w}{i} - 1 \right)
\end{equation*}
and set $\lambda$ equal to the right hand side of this inequality. With $\lambda$ in hand, the Gibbs condition implies that $\VE{x}{i} = \max \{\VE{y}{i} - \lambda, 0\}$.
It follows that projection onto $U$ can be accomplished in $O(d \log d)$ operations dominated by sorting. Algorithm~\ref{alg:project_simplex} displays pseudocode for projection onto
$U$.

\begin{algorithm}[t]
	\begin{algorithmic}[0]
		\State $\V{w} \leftarrow \textsc{sort\_descending}(\V{y}).$
		\State $k \leftarrow \max \left \{ j  : \VE{w}{j}  >  \frac{1}{j} \left( \sum_{i=1}^j \VE{w}{i} \right)\right \}$
		\State $\lambda \leftarrow \frac{1}{k} \left( \sum_{i=1}^{k} \VE{w}{i} \right)$
		\State $\VE{x}{i} \leftarrow \max\{\VE{y}{i} - \lambda, 0\}$.
	\end{algorithmic}
	\caption{\, \textsc{Projection onto simplex}}	
	 \label{alg:project_simplex}
\end{algorithm}

Armed with these results, we now describe how to solve a continuous relaxation of Sudoku by the method of alternating projections.   In the relaxed version of the problem, we imagine generating candidate solutions by random sampling.  Each cell $(i,j)$ is assigned a sampling distribution $\TE{p}{ijk}= \Pr(S_{ij}=k)$ for choosing a random deviate $S_{ij} \in \{1,\ldots,9\}$ to populate the cell. If a numerical clue $k$ occupies cell $(i,j)$, then we set $p_{ijl}=1$ for $l=k$ and 0 otherwise. A matrix of sampled deviates $\M{S}$ constitutes a candidate solution. It seems reasonable to demand that the average puzzle obey the constraints.  Once we find a feasible 3-dimensional tensor $\T{P} = (\TE{p}{ijk})$ obeying the constraints, a good heuristic for generating an integer solution $\Mhat{S}$ is to put
\begin{eqnarray*}
\hat{s}_{ij}  & = & \underset{k \in \{ 1, \ldots, 9 \}}{\max} \TE{P}{ijk}.
\end{eqnarray*}
In other words, we impute the most probable integer to each unknown cell $(i,j)$. It is easy to construct counterexamples where imputation of the most probable integer from a feasible tensor $\T{P}$ of the relaxed problem fails to solve the Sudoku puzzle.

In any case, the remaining agenda is to specify the constraints and the corresponding projection operators.  The requirement that each digit appear in each row on average once amounts to the constraint $\sum_{j=1}^9 \TE{p}{ijk} = 1$ for all $i$ and $k$ between 1 and 9.  There are 81 such constraints. The requirement that each digit appear in each column on average once amounts to the constraint $\sum_{i=1}^9 \TE{p}{ijk} = 1$ for all $j$ and $k$ between 1 and 9. Again, there are 81 such constraints. The requirement that each digit appear in each subgrid on average once amounts to the constraint $\sum_{j=1}^3 \sum_{j=1}^3 \TE{p}{a+i,b+j,k} = 1$ for all $k$ between 1 and 9 and all $a$ and $b$ chosen from the set $\{0,3,6\}$.  This contributes another 81 constraints. Finally, the probability constraints $\sum_{k=1}^9 \TE{p}{ijk} = 1$ for all $i$ and $j$ between 1 and 9 contribute 81 more affine constraints.  Hence, there are a total of $324$ affine constraints on the $9^3 = 729$ parameters.  In addition there are 729 nonnegativity constraints $\TE{p}{ijk} \ge 0$.

Every numerical clue voids several constraints. For example, if the digit $7$ is mandated for cell (9,2), then we must take 
$\TE{p}{927} = 1$, $\TE{p}{92k} = 0$ for $k \ne 7$, $\TE{p}{i27} = 0$ for all $i \not = 9$,  $\TE{p}{9j7} = 0$ for all $j \not = 2$, and $\TE{p}{ij7} = 0$ for all other pairs $(i,j)$ in the (3,1) subgrid. In carrying out alternating projection, we eliminate the corresponding variables.  With this proviso, we cycle through the simplex projections 
summarized in Algorithm~\ref{alg:project_simplex}.
The process is very efficient but slightly tedious to code.  For the sake of brevity we omit the remaining details. All code used to generate the subsequent 
results are available at \url{https://github.com/echi/Sudoku}, and we direct the interested reader there.

\section{Comparisons}
\label{sec:experiments}

We generated test puzzles from code available online \cite{Wan2007} and discarded puzzles that could be completely solved by filling in entries directly implied by the initial clues. This left 87 easy puzzles, 130 medium puzzles, and 100 hard puzzles.  We also downloaded an additional 95 very hard benchmark puzzles \cite{Coc2012,Ste2005}. In simulated annealing, the temperature $\tau$ was initialized to 200 and lowered by a factor of 0.99 after every 50 steps. 
We allowed at most $2 \times 10^{5}$ iterations and reset the temperature to 200 if a solution had not been found after $10^{5}$ iterations.  For the alternating projection algorithm, we commenced projecting from the origin $\bf 0$. 

\begin{table}[th]
  \centering
  \begin{tabular}{lcccc}
  \hline \hline
	& Alt. Projection & Sim. Annealing & Backtracking & Number of Puzzles\\ \hline 
	Easy & 0.85 & 1.00 & 1.00 & 87 \\
	Medium & 0.89 & 1.00 & 1.00 & 130 \\
	Hard & 0.72 & 0.97 & 1.00 & 100 \\
	Top 95 & 0.41 & 0.03 & 1.00 & 95 \\ \hline
  \end{tabular}
  \caption{Success rates for solving puzzles of varying difficulty.}
  \label{tab:comparisons}
\end{table}

\begin{table}[th]
\centering
\begin{tabular}{clccc}
\toprule
\midrule
& & \multicolumn{3}{c}{CPU Time (sec)}  \\
\cmidrule(r){3-5}
& & Alt. Projection & Sim. Annealing & Backtracking \\ 
\midrule
\multicolumn{1}{c}{\multirow{4}{*}{Easy}} &
\multicolumn{1}{l}{Minimum} & 0.032 & 0.006 & 0.007 \\ 
\multicolumn{1}{l}{}                        &
\multicolumn{1}{l}{Median} & 0.041 & 0.021 & 0.008 \\ 
\multicolumn{1}{l}{}                        &
\multicolumn{1}{l}{Mean} & 0.052 & 0.112 & 0.008 \\
\multicolumn{1}{l}{}                        &
\multicolumn{1}{l}{Maximum} & 0.237 & 0.970 & 0.009 \\ \\
\multicolumn{1}{c}{\multirow{4}{*}{Medium}} &
\multicolumn{1}{l}{Minimum} & 0.032 & 0.007 & 0.007 \\  
\multicolumn{1}{l}{}                        &
\multicolumn{1}{l}{Median} & 0.051 & 0.037 & 0.008 \\
\multicolumn{1}{l}{}                        &
\multicolumn{1}{l}{Mean} & 0.062 & 0.231 & 0.008 \\ 
\multicolumn{1}{l}{}                        &
\multicolumn{1}{l}{Maximum} & 0.269 & 3.36 & 0.010 \\ \\
\multicolumn{1}{c}{\multirow{4}{*}{Hard}} &
\multicolumn{1}{l}{Minimum} & 0.033 & 0.008 & 0.008 \\ 
\multicolumn{1}{l}{}                        &
\multicolumn{1}{l}{Median} & 0.110 & 0.753 & 0.008 \\ 
\multicolumn{1}{l}{}                        &
\multicolumn{1}{l}{Mean} & 0.159 & 1.104 & 0.009 \\ 
\multicolumn{1}{l}{}                        &
\multicolumn{1}{l}{Maximum} & 0.525 & 7.204 & 0.031 \\ 
\bottomrule
\end{tabular}
  \caption{Summary statistics on the run times for different methods on puzzles of varying difficulty. For the alternating projection and simulated annealing techniques, only successfully solved puzzles are included in the statistics.}
  \label{tab:run_times}
\end{table}

Backtracking successfully solved all puzzles. Table~\ref{tab:comparisons} shows the fraction of puzzles the two heuristics were able to successfully complete. Table~\ref{tab:run_times} records summary statistics for the CPU time taken by each method for each puzzle category. All computations were done on an iMac computer with a 3.4 GHz Intel Core i7 processor and 8 GB of RAM. We implemented the alternating projection and simulated annealing algorithms in Fortran 95. For backtracking we relied on the existing implementation in C. 

The comparisons show that backtracking performs best, and for the vast majority of $9 \times 9$ Sudoku problems it is probably going to be hard to beat. Simulated annealing 
finds the solution except for a handful of the most challenging paper and pencil problems, but its maximum run times are unimpressive. While alternating projection does not perform as well on the pencil and paper problems compared to the other two algorithms, it does not do terribly either. Moreover, we see hints of the tables turning on the hard puzzles.

\begin{figure}[t]
\centering
\begin{tikzpicture}[scale=0.75]
    \draw[black!50] (0,0) grid (9,9);
	\draw[black!100] (0,0) rectangle (3,3);
	\draw[black!100] (3,0) rectangle (6,3);
	\draw[black!100] (6,0) rectangle (9,3);
	\draw[black!100] (0,3) rectangle (3,6);
	\draw[black!100] (3,3) rectangle (6,6);
	\draw[black!100] (6,3) rectangle (9,6);
	\draw[black!100] (0,6) rectangle (3,9);
	\draw[black!100] (3,6) rectangle (6,9);
	\draw[black!100] (6,6) rectangle (9,9);	
\node[regular polygon, regular polygon sides=4,fill=gray!20,scale=0.75] at (0.5,8.5) {$4$};
\node[regular polygon, regular polygon sides=4,scale=0.75] at (1.5,8.5) {$1$};
\node[regular polygon, regular polygon sides=4,scale=0.75] at (2.5,8.5) {$7$};
\node[regular polygon, regular polygon sides=4,scale=0.75] at (3.5,8.5) {$9$};
\node[regular polygon, regular polygon sides=4,scale=0.75] at (4.5,8.5) {$6$};
\node[regular polygon, regular polygon sides=4,fill=gray!70,scale=0.75] at (5.5,8.5) {$2$};
\node[regular polygon, regular polygon sides=4,fill=gray!20,scale=0.75] at (6.5,8.5) {$8$};
\node[regular polygon, regular polygon sides=4,scale=0.75] at (7.5,8.5) {$3$};
\node[regular polygon, regular polygon sides=4,fill=gray!20,scale=0.75] at (8.5,8.5) {$5$};
\node[regular polygon, regular polygon sides=4,scale=0.75] at (0.5,7.5) {$6$};
\node[regular polygon, regular polygon sides=4,fill=gray!20,scale=0.75] at (1.5,7.5) {$3$};
\node[regular polygon, regular polygon sides=4,scale=0.75] at (2.5,7.5) {$2$};
\node[regular polygon, regular polygon sides=4,scale=0.75] at (3.5,7.5) {$1$};
\node[regular polygon, regular polygon sides=4,fill=gray!70,scale=0.75] at (4.5,7.5) {$5$};
\node[regular polygon, regular polygon sides=4,scale=0.75] at (5.5,7.5) {$8$};
\node[regular polygon, regular polygon sides=4,fill=gray!20,scale=0.75] at (6.5,7.5) {$7$};
\node[regular polygon, regular polygon sides=4,scale=0.75] at (7.5,7.5) {$4$};
\node[regular polygon, regular polygon sides=4,scale=0.75] at (8.5,7.5) {$9$};
\node[regular polygon, regular polygon sides=4,scale=0.75] at (0.5,6.5) {$9$};
\node[regular polygon, regular polygon sides=4,scale=0.75] at (1.5,6.5) {$5$};
\node[regular polygon, regular polygon sides=4,scale=0.75] at (2.5,6.5) {$8$};
\node[regular polygon, regular polygon sides=4,scale=0.75] at (3.5,6.5) {$7$};
\node[regular polygon, regular polygon sides=4,scale=0.75] at (4.5,6.5) {$3$};
\node[regular polygon, regular polygon sides=4,scale=0.75] at (5.5,6.5) {$4$};
\node[regular polygon, regular polygon sides=4,scale=0.75] at (6.5,6.5) {$6$};
\node[regular polygon, regular polygon sides=4,scale=0.75] at (7.5,6.5) {$1$};
\node[regular polygon, regular polygon sides=4,scale=0.75] at (8.5,6.5) {$2$};
\node[regular polygon, regular polygon sides=4,scale=0.75] at (0.5,5.5) {$8$};
\node[regular polygon, regular polygon sides=4,fill=gray!20,scale=0.75] at (1.5,5.5) {$2$};
\node[regular polygon, regular polygon sides=4,scale=0.75] at (2.5,5.5) {$5$};
\node[regular polygon, regular polygon sides=4,scale=0.75] at (3.5,5.5) {$4$};
\node[regular polygon, regular polygon sides=4,scale=0.75] at (4.5,5.5) {$9$};
\node[regular polygon, regular polygon sides=4,scale=0.75] at (5.5,5.5) {$7$};
\node[regular polygon, regular polygon sides=4,scale=0.75] at (6.5,5.5) {$3$};
\node[regular polygon, regular polygon sides=4,fill=gray!20,scale=0.75] at (7.5,5.5) {$6$};
\node[regular polygon, regular polygon sides=4,scale=0.75] at (8.5,5.5) {$1$};
\node[regular polygon, regular polygon sides=4,scale=0.75] at (0.5,4.5) {$3$};
\node[regular polygon, regular polygon sides=4,scale=0.75] at (1.5,4.5) {$9$};
\node[regular polygon, regular polygon sides=4,scale=0.75] at (2.5,4.5) {$1$};
\node[regular polygon, regular polygon sides=4,scale=0.75] at (3.5,4.5) {$5$};
\node[regular polygon, regular polygon sides=4,fill=gray!20,scale=0.75] at (4.5,4.5) {$8$};
\node[regular polygon, regular polygon sides=4,scale=0.75] at (5.5,4.5) {$6$};
\node[regular polygon, regular polygon sides=4,fill=gray!20,scale=0.75] at (6.5,4.5) {$4$};
\node[regular polygon, regular polygon sides=4,scale=0.75] at (7.5,4.5) {$2$};
\node[regular polygon, regular polygon sides=4,scale=0.75] at (8.5,4.5) {$7$};
\node[regular polygon, regular polygon sides=4,scale=0.75] at (0.5,3.5) {$7$};
\node[regular polygon, regular polygon sides=4,fill=gray!20,scale=0.75] at (1.5,3.5) {$4$};
\node[regular polygon, regular polygon sides=4,scale=0.75] at (2.5,3.5) {$6$};
\node[regular polygon, regular polygon sides=4,scale=0.75] at (3.5,3.5) {$3$};
\node[regular polygon, regular polygon sides=4,fill=gray!20,scale=0.75] at (4.5,3.5) {$1$};
\node[regular polygon, regular polygon sides=4,fill=gray!70,scale=0.75] at (5.5,3.5) {$2$};
\node[regular polygon, regular polygon sides=4,scale=0.75] at (6.5,3.5) {$5$};
\node[regular polygon, regular polygon sides=4,scale=0.75] at (7.5,3.5) {$9$};
\node[regular polygon, regular polygon sides=4,scale=0.75] at (8.5,3.5) {$8$};
\node[regular polygon, regular polygon sides=4,scale=0.75] at (0.5,2.5) {$2$};
\node[regular polygon, regular polygon sides=4,scale=0.75] at (1.5,2.5) {$8$};
\node[regular polygon, regular polygon sides=4,scale=0.75] at (2.5,2.5) {$9$};
\node[regular polygon, regular polygon sides=4,fill=gray!20,scale=0.75] at (3.5,2.5) {$6$};
\node[regular polygon, regular polygon sides=4,fill=gray!70,scale=0.75] at (4.5,2.5) {$5$};
\node[regular polygon, regular polygon sides=4,fill=gray!20,scale=0.75] at (5.5,2.5) {$3$};
\node[regular polygon, regular polygon sides=4,scale=0.75] at (6.5,2.5) {$1$};
\node[regular polygon, regular polygon sides=4,fill=gray!20,scale=0.75] at (7.5,2.5) {$7$};
\node[regular polygon, regular polygon sides=4,scale=0.75] at (8.5,2.5) {$4$};
\node[regular polygon, regular polygon sides=4,fill=gray!20,scale=0.75] at (0.5,1.5) {$5$};
\node[regular polygon, regular polygon sides=4,scale=0.75] at (1.5,1.5) {$7$};
\node[regular polygon, regular polygon sides=4,fill=gray!20,scale=0.75] at (2.5,1.5) {$3$};
\node[regular polygon, regular polygon sides=4,fill=gray!20,scale=0.75] at (3.5,1.5) {$2$};
\node[regular polygon, regular polygon sides=4,scale=0.75] at (4.5,1.5) {$4$};
\node[regular polygon, regular polygon sides=4,fill=gray!20,scale=0.75] at (5.5,1.5) {$1$};
\node[regular polygon, regular polygon sides=4,scale=0.75] at (6.5,1.5) {$9$};
\node[regular polygon, regular polygon sides=4,scale=0.75] at (7.5,1.5) {$8$};
\node[regular polygon, regular polygon sides=4,scale=0.75] at (8.5,1.5) {$6$};
\node[regular polygon, regular polygon sides=4,fill=gray!20,scale=0.75] at (0.5,0.5) {$1$};
\node[regular polygon, regular polygon sides=4,scale=0.75] at (1.5,0.5) {$6$};
\node[regular polygon, regular polygon sides=4,fill=gray!20,scale=0.75] at (2.5,0.5) {$4$};
\node[regular polygon, regular polygon sides=4,scale=0.75] at (3.5,0.5) {$8$};
\node[regular polygon, regular polygon sides=4,scale=0.75] at (4.5,0.5) {$7$};
\node[regular polygon, regular polygon sides=4,scale=0.75] at (5.5,0.5) {$9$};
\node[regular polygon, regular polygon sides=4,scale=0.75] at (6.5,0.5) {$2$};
\node[regular polygon, regular polygon sides=4,scale=0.75] at (7.5,0.5) {$5$};
\node[regular polygon, regular polygon sides=4,scale=0.75] at (8.5,0.5) {$3$};
\end{tikzpicture}
\caption{A typical local minimum that traps simulated annealing in a top 95 puzzle. Clues are shaded light gray. There are two
column constraint violations caused by the cells shaded dark gray. The local minimum is deep in the sense that all one-step swaps  
result in further constraint violations. \label{fig:top95_1}}
\end{figure}

\begin{figure}
\centering
\includegraphics[scale=0.75]{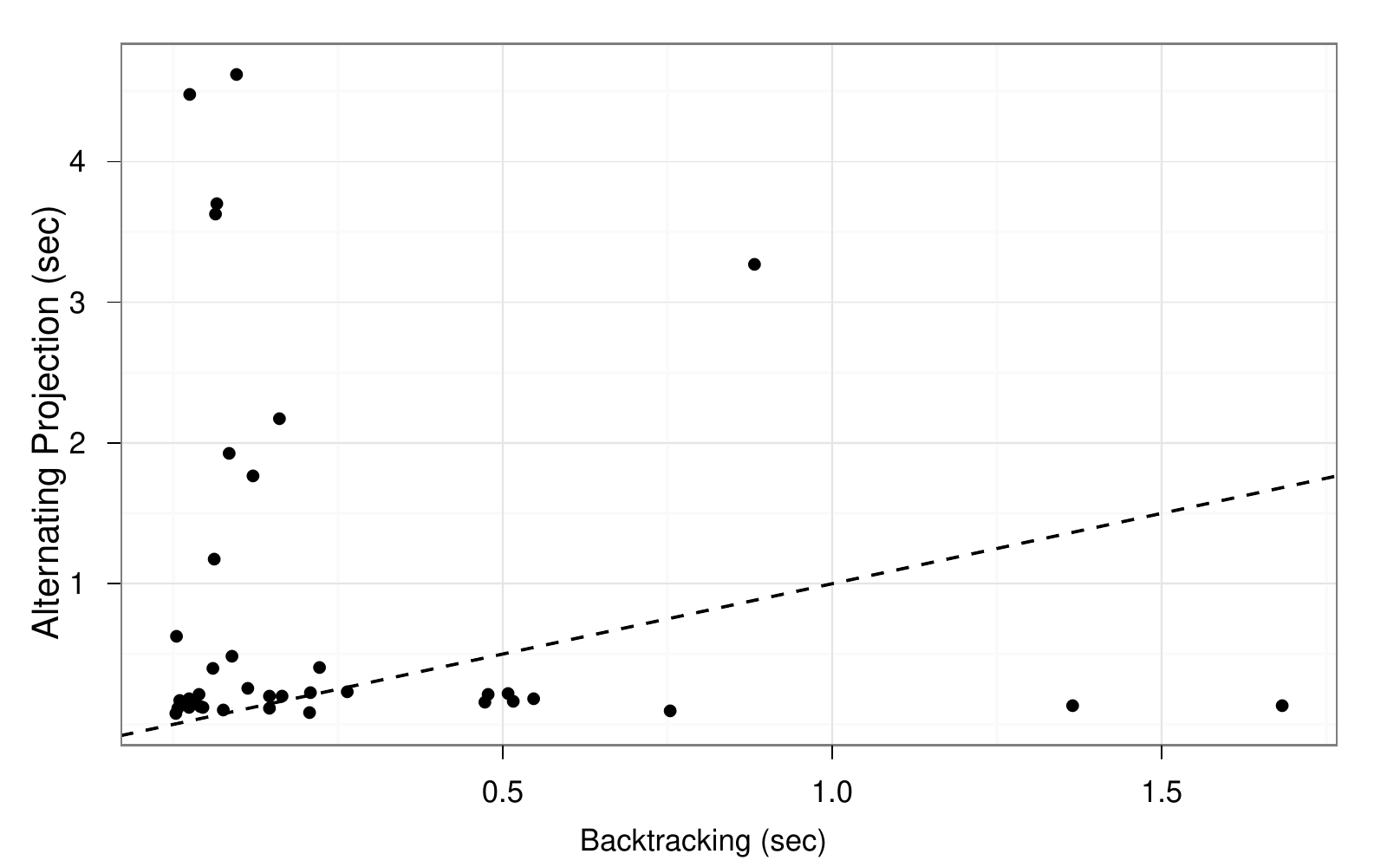}
\caption{Scatterplot of solution times for the top 95 benchmark problems.}
\label{fig:top}
\end{figure}

Simulated annealing struggles mightily on the 95 benchmark puzzles. Closer inspection of individual puzzles reveals that these very hard puzzles admit many local minima with just a few constraint violations. Figure~\ref{fig:top95_1} shows a typical local minimum that traps the simulated annealing algorithm. Additionally, something curious is happening in Figure~\ref{fig:top}, which plots CPU solution times for alternating projection versus backtracking. Points below the dashed line indicate puzzles that the method of alternating projection solves more efficiently than backtracking. It appears that when the method of alternating projections finds correct solutions to very hard problems, it tends to find them more quickly than backtracking.

\section{Discussion}
\label{sec:discussion} 

It goes almost without saying that students of the mathematical sciences should be taught a variety of solution techniques for combinatorial problems. Because Sudoku puzzles are easy to state and culturally neutral, they furnish a good starting point for the educational task ahead.  It is important to stress the contrast between exact strategies that scale poorly with problem size and approximate strategies that adapt more gracefully.  The performance of the alternating projection algorithm on the benchmark tests suggest it may have a role in solving much harder combinatorial problems.  Certainly, the electrical engineering community takes this attitude, given the close kinship of Sudoku puzzles to problems in coding theory \cite{ErlChaGor2009, GunMoo2012, MooGun2006, MooGunKup2009}.

One can argue that algorithm development has assumed a dominant role within the mathematical sciences.   Three inter-related trends are feeding this phenomenon. First, computing power continues to grow.  Execution times are dropping, and computer memory is getting cheaper.  Second, good computing simply tempts scientists to tackle larger data sets.  Third, certain fields, notably communications, imaging, genomics, and economics generate enormous amounts of data.  All of these fields create problems in combinatorial optimization. For instance, modern DNA sequencing is still challenged by the phase problem of discerning the maternal and paternal origin of genetic variants.  Computation is also being adopted more often as a means of proving propositions. The claim that at least 17 numerical clues are needed to ensure uniqueness of a Sudoku solution has apparently been proved using intelligent brute force \cite{McGTugCiv2012}.  Mathematical scientists need to be aware of computational developments outside their classical application silos.  Importing algorithms from outside fields is one of the quickest means of refreshing an existing field.

\section*{Acknowledgments}
We thank Peter Cock for helpful comments and pointing us to the source and author of the ``Top 95" puzzles.

\bibliographystyle{siammod}
\bibliography{sudoku}

\end{document}